\documentclass[10pt]{article}
\usepackage{graphicx}
\usepackage{biograph}        
\usepackage{amsmath,amsfonts,amssymb,amsthm}
\usepackage{latexsym,xypic,eucal,mathrsfs,setspace}
\newcommand{\n}{\mathfrak{N}}
\newcommand{\g}{\mathfrak{g}}

\newcommand{\cart}{\mathfrak{h}}

\newcommand{\goo}{\mathfrak{g}_{{0}}}
\newcommand{\goc}{\g_{{0}}^{{\mathbb{C}}}}
\newcommand{\ngam}{\Gamma \backslash N}

\newcommand{\linl}{\lambda\in\Lambda}
\newcommand{\aone}{\alpha_1}
\newcommand{\atwo}{\alpha_2}
\newcommand{\ai}{\alpha_i}
\newcommand{\an}{\alpha_n}
\newcommand{\anone}{\alpha_{n-1}}
\newcommand{\athree}{\alpha_3}
\newcommand{\afour}{\alpha_4}
\newcommand{\afive}{\alpha_5}
\newcommand{\asix}{\alpha_6}
\newcommand{\aseven}{\alpha_7}
\newcommand{\aeight}{\alpha_8}
\newcommand{\mone}{m_1}
\newcommand{\mtwo}{m_2}
\newcommand{\mn}{m_n}

\newcommand{\mthree}{m_3}
\newcommand{\mfour}{m_4}
\newcommand{\mfive}{m_5}
\newcommand{\msix}{m_6}
\newcommand{\mseven}{m_7}
\newcommand{\meight}{m_8}
\newcommand{\mi}{m_i}
\newcommand{\mimone}{m_{i-1}}
\newcommand{\mipone}{m_{i+1}}

\newcommand{\klam}{K_{\lambda,\mu}}
\newcommand{\klami}{K_{\lambda,\mu_i}}
\newcommand{\kzeta}{K_{\zeta,\mu}}

\newtheorem{theorem}{Theorem}[section]
\newtheorem{lemma}[theorem]{Lemma}

\newtheorem{proposition}[theorem]{Proposition}

\newtheorem{definition}[theorem]{Definition}
\newtheorem{sublemma}{Lemma}[subsection]

\begin{document}

\title{The multiplicity of weights in nonprimitive pairs of weights
}


\author{Rachelle C. DeCoste     \\
Department of Mathematics\\
 Wheaton College, Norton, MA 02766 \\
              Tel.: 508-286-5842\\
             decoste\underline{ }rachelle@wheatoncollege.edu             
}



\date{\today}

\maketitle

\begin{abstract}
For each type of classical Lie algebra, we list the dominant highest weights $\zeta$ for which $(\zeta;\mu_i)$ is not a primitive pair and the weight space $V_{\mu_i}$ has dimension one where $\mu_i$ are the highest long and short roots in each case.  These dimension one weight spaces lead to examples of nilmanifolds for which we cannot prove or disprove the density of closed geodesics.
\end{abstract}

\section{Introduction}
\label{intro}


In our study \cite{decoste} of the distribution of closed geodesics on nilmanifolds, we considered manifolds arising from a Lie group $N$ with an associated Lie algebra $\n$ constructed from an irreducible representation of a compact semisimple Lie algebra $\goo$ on a real finite dimensional vector space $U$.  The nilmanifolds considered, $\ngam$, are those such that $\Gamma$ arises from a Chevalley rational structure on $\n$.   The main result of that study classified such nilmanifolds as having the density of closed geodesics property if all roots of $\g=\goc$ were weights of $V=U^\mathbb{C}$ with multiplicity greater than or equal to two.



In \cite{decoste} we reduce the multiplicity question to $\g$ simple, thus throughout this article, we will assume that $\g$ is a complex simple Lie algebra with a fixed base $\Delta$ of positive simple roots determined by a Cartan subalgebra $\cart$.  Let $V=V(\lambda)$ denote a finite dimensional irreducible $g-$module with highest weight $\lambda$.  The multiplicity of a weight $\mu$ is defined to be the dimension of the weight space $V_{\mu}\subseteq V$, and is denoted $K_{\lambda,\mu}$.  By standard results of Lie theory (cf \cite{humph}),  each root of $\g$ is conjugate to the highest short root $\mu_1$ or the highest long root $\mu_2$, thus we can consider only these roots when finding the dimension of the weight spaces of interest.  The results of this paper provide for each classical Lie algebra type, all highest weights for which the highest short or long roots give rise to weight spaces of dimension one when the root and weight are nonprimitive pairs as defined below.  The case of primitive pairs is completely answered in \cite{decoste}.  Thus this paper provides many cases in \cite{decoste} for which the density of closed geodesics cannot be shown with the traditional methods used.  These exceptional cases provide examples for which the distribution of closed geodesics is unknown and still being investigated.  They may provide unique examples of nilmanifolds satisfying necessary conditions, but not having the density of closed geodesics.

\begin{definition} For $\g$ simple, define a pair $(\lambda; \mu)$ of weights in $\Lambda^+$ to be \emph{primitive} if $(\lambda-\mu)$ written as the sum of simple roots has all
positive integer coefficients.  
\end{definition}

Thus the pair $(\lambda;\mu)$ is said to be \emph{nonprimitive} if in the sum $\lambda-\mu$, at least one simple root has a zero coefficient.  The weights for which the highest roots give rise to weight spaces of dimension one in the primitive pair case as discussed in \cite{decoste}  were found using Theorem \ref{bzthm} below which is also applied repeatedly in this article.

By \cite{bz} we will be able to reduce to the primitive case to find all weight spaces of dimension one.  Thus the following result will be the basis of our determination of all dominant weights $\lambda$ such that $K_{\lambda,\mu_i}=1$ for $i=1,2$.  
In the notation of \cite{bz}, $\mathbb{Z}_+$ is the set of all nonnegative integers and $\{\omega_i\}$ is the set of fundamental dominant weights relative to $\Delta$ (found in Table 1, page 69 of \cite{humph}).  Additionally, the partial ordering of weights $\lambda\succ\mu$ means that $\lambda-\mu$ is a linear combination of simple roots with nonnegative coefficients.

\begin{theorem} [\cite{bz}, Theorem 1.3]  All primitive pairs $(\lambda; \mu)$ such that $K_{\lambda,\mu}=1$, up to
isomorphism of Dynkin diagrams, are exhausted by the following list: 
\begin{enumerate}
\item $A_n$ $(n\geq 1)$: $\lambda=l\omega_1$, $\displaystyle \mu=\sum_{1\leq i \leq n}a_i\omega_i$ 
  where $a_i\in\mathbb{Z}_+$ and \\ 
$\displaystyle (l-\sum_{1\leq i\leq n} ia_i)\in (n+1)\mathbb{N}$ 
\item $B_n$ $(n\geq 2)$: $\lambda=l\omega_1$, $\displaystyle \mu=\sum_{1\leq i \leq n}a_i\omega_i$ 
 where $a_i\in\mathbb{Z}_+$
is even and \\
$\displaystyle (l-1)=\sum_{1\leq i \leq n-1}ia_i+\frac{na_n}{2}$
\item $G_2$: $\lambda=l\omega_2$, $\displaystyle \mu=a_1\omega_1+a_2\omega_2$
  where $a_1,\ a_2 \in \mathbb{Z}_+$,
and $3l-1=2a_1+3a_2$
\item $G_2$:   $\lambda= \omega_1$, $\mu=0$.
\end{enumerate}
\label{bzthm}
\end{theorem}

For each $\g$ of classical type, we first identify those highest dominant weights $\lambda$ having the property that $\lambda\neq\mu_i$ and $(\lambda; \mu_i)$ is not a primitive pair for $\mu_i$, the highest short or long root $i=1,2$ respectively. Once we have identified such $\lambda$, we reduce to a primitive pair by another result of \cite{bz} and then use  Theorem \ref{bzthm} to determine whether $\klami=1$.  For review, we list the highest short and long roots in Table \ref{domwts}\cite{humph}.

\begin{table}[h]
\caption{Highest short and long roots}
\centering
\label{domwts}
\begin{tabular}{|c|c|} \hline
Classical Lie algebra type &  \\ \hline
$A_n$ & $\alpha_1+\alpha_2+\cdots+\alpha_n$\\ \hline
$B_n$ & $\alpha_1+\alpha_2+\cdots +\alpha_n$ \\ 
 	& $\alpha_1+2\alpha_2\cdots +2\alpha_n$ \\ \hline
$C_n$ & $\alpha_1+2\alpha_2+2\alpha_3+\cdots+2\alpha_{n-1}+\alpha_n$\\
	& $2\alpha_1+2\alpha_2+2\alpha_3+\cdots+2\alpha_{n-1}+\alpha_n$\\ \hline
$D_n$ & $\alpha_1+2\alpha_2+\cdots +2\alpha_{n-2}+\alpha_{n-1}+\alpha_n$\\ \hline
$E_6$	& $\alpha_1+2\alpha_2+2\alpha_3+3\alpha_4+2\alpha_5+\alpha_6$\\ \hline
$E_7$	& $2\alpha_1+2\alpha_2+3\alpha_3+4\alpha_4+3\alpha_5+2\alpha_6+\alpha_7$\\ \hline
$E_8$	& $2\alpha_1+3\alpha_2+4\alpha_3+6\alpha_4+5\alpha_5+4\alpha_6+3\alpha_7+2\alpha_8$\\ \hline
$F_4$ & $\alpha_1+2\alpha_2+3\alpha_3+2\alpha_4$ \\ 
	&   $2\alpha_1+3\alpha_2+4\alpha_3+2\alpha_4$ \\ \hline
$G_2$ & $2\alpha_1+\alpha_2$ \\ 
	& $3\alpha_1+2\alpha_2$ \\ \hline
\end{tabular}
\end{table}

The remainder of this paper is the proof of the following result, considering case by case each class of simple Lie algebra.

\begin{theorem}
Consider $\mu=\mu_i$, $i=1,2$, the highest short and long roots of $\g$.  All nonprimitive pairs $(\zeta;\mu_i)$ such that $K_{\zeta,\mu_i}=1$, up to isomorphism of Dynkin diagrams are exhausted by the Table \ref{lprime}.\begin{table}[h]
\caption{Theorem \ref{main}}
\centering
\label{lprime}
\begin{tabular}{|c|c|} \hline
Lie algebra type & dominant highest weight $\zeta$ \\ \hline
$A_n$  & $\alpha_1+2\alpha_2$, $2\alpha_1+\alpha_2$ $(n=2)$  \\
	& $\alpha_1+2\alpha_2+\alpha_3$ $(n=3)$  \\ \hline
$B_n$  $(n\geq 2)$ & $\alpha_1+2\alpha_2+m_3\alpha_3$, $m_3\geq 3$ $(n=3)$\\
	& $2\alpha_1+2\alpha_2\cdots +2\alpha_n$\\ 
	& $\alpha_1+2\alpha_2+3\alpha_3+3\alpha_4+\cdots +3\alpha_n$ $(n\geq4)$\\ \hline
$D_n$ $(n\geq 4)$  &$2\alpha_1+2\alpha_2+\cdots+2\alpha_{n-2}+\alpha_{n-1}+\alpha_n$\\ 
&  $\alpha_1+2\alpha_2+2\alpha_3+\alpha_4$, $\alpha_1+2\alpha_2+\alpha_3+2\alpha_4$ $(n=4)$\\ \hline
$G_2$ & $4\alpha_1+2\alpha_2$ \\ \hline
\end{tabular}
\end{table}
\label{main}
\end{theorem}

Before continuing, a few more definitions and the following result are necessary.  For any subset $S$ of the set of simple roots $\Delta$, define $\g(S)$ to be the subalgebra of $\g$ generated by the root subspaces $\g_\beta$ and $\g_{-\beta}$ for all $\beta\in S$.  We define the projection map $p=p_S$ to be the natural projection of the set of weights of $\g$ to the set of weights of $\g(S)$.  The Lie algebra $\g(S)$ is known to be semisimple.

\begin{proposition} \cite[Proposition 2.4]{bz}  \begin{enumerate}
\item Let $S$ be a subset of simple roots.  Let $\linl$ be an element such that the expansion of the weight $(\lambda-\mu)$ in terms of simple roots involves only elements of $S$.  Then $\klam=K_{p(\lambda),p(\mu)}$.
\item Under the assumptions of part 1, let $S_1,\dots,S_k$ be all the connected components of the set $S$ in the Dynkin diagram of the system of positive roots, and let $\lambda_i=p_{S_i}(\lambda)$ and $\mu_i=p_{S_i}(\mu)$.  Then $\klam=\Pi_{1\leq i\leq k}K_{\lambda_i,\mu_i}$.
\end{enumerate}
\label{bztwo}
\end{proposition}

In our results, since we have reduced to considering $\g$ simple, let $S$ be the set of simple roots that occur with nonzero coefficients in the difference $\lambda-\mu$ and let $p=p_S$.  We can then apply Proposition \ref{bztwo} and Theorem \ref{bzthm} to determine if $K_{p(\lambda),p(\mu)}=1$.  Thus for each Lie algebra and each highest weight $\lambda$ where $(\lambda;\mu)$ is not a primitive pair, we must identify $S$ and then determine $\g(S)$.  To find $\g(S)$, we can simply consider the Dynkin diagram of the root system of $S$.

In many cases $S=\{\alpha_i\}$ and then $\g(S)\cong A_1$.  Thus the following lemma will be used frequently as we continue our discussion of nonprimitive pairs.

\begin{lemma}  Let $\g$ be of type $A_1$ with root $\alpha_1$.  Then $K_{2\alpha_1,\alpha_1}=1$.\label{klem}\end{lemma}

\begin{proof} Let $\lambda=2\alpha_1$ and $\mu=\alpha_1$.  By Theorem \ref{bzthm} for the case $A_1$, $\klam=1$ if $\lambda=l\omega_1$, $\mu=a_1\omega_1$ with $a_1\in\mathbb{Z}_+$ and $l-a_1\in 2\mathbb{N}$.  Since $\mu=\alpha_1=2\omega_1$, $a_1=2$ and since $\lambda=2\alpha_1=4\omega_1$, $l=4$.  Clearly then $l-a_1\in 2\mathbb{N}$ and $\klam=1$. \end{proof}

In each case of classical Lie algebra, we assume the zero weight space of our representation $V$ is nontrivial, as necessary in \cite{decoste}.  Therefore every dominant weight $\lambda$ is of the form $\lambda=\sum \mi\ai$ where $\mi$ must satisfy inequalities that arise from the conditions $\langle \lambda, \ai \rangle \geq 0$ for $i=1,\dots,n$.  Also $\mi>0$ for all $i$ by the following lemma of \cite{decoste}.

\begin{lemma}  Let $\mu\in\Lambda^+$ and suppose that $\displaystyle\mu=\sum_{i=1}^nm_i\alpha_i$ for integers $\{m_i\}$.  Then $m_k>0$ for all $k$.
\label{posint}
\end{lemma}

This lemma is used in the proof of the following necessary proposition, found in \cite{decoste}.

\begin{proposition} Let $V$ be an irreducible $\g-$module with nontrivial zero weight space $V_0$.  Let $\lambda\in\Lambda(V)^+$ be the highest weight.  Then
\begin{enumerate}
\item $\displaystyle\lambda=\sum_{i=1}^n p_i\alpha_i$ for suitable positive integers $p_i$.
\item If $\mu\in\Lambda(V)$, then $\displaystyle\mu=\sum_{i=1}^n m_i\alpha_i$, $\mi\in\mathbb{Z}$.  Furthermore, if $\mu\in\Lambda^+(V)$, then the integers $\{m_i\}$ are all positive.
\item At least one root of $\g$ is a weight.
\end{enumerate}
\label{atleast1} \end{proposition}

\section{Nonprimitive pairs for $A_n$}

In this case, recall that $\mu_1=\mu_2=\mu=\alpha_1+\cdots+\alpha_n$ is the dominant weight that is conjugate to all roots.  We will find all dominant weights $\zeta$ such that $\zeta\succ\mu$ and $(\zeta;\mu)$ is not a primitive pair.  Next we determine for which of these weights $K_{\zeta,\mu}=1$.  

\begin{lemma} For $\g$ of type $A_n$ with $n\geq 4$, there are no highest weights $\zeta$ such that $\zeta\succ\mu$, $(\zeta;\mu)$ is not a primitive pair and $K_{\zeta,\mu}=1$.  For $n=2,3$ the following exceptional cases occur such that $K_{\zeta_i,\mu}=1$.
\begin{enumerate}
\item $n=2 \  \ \ \zeta_1=\alpha_1+2\alpha_2$, $\zeta_2=2\alpha_1+\alpha_2$
\item $n=3 \  \ \ \zeta_3=\alpha_1+2\alpha_2+\alpha_3$
\end{enumerate}
\label{annonprim}
\end{lemma}

\begin{lemma} If $S$ is a set of $k$ consecutive simple roots in $A_n$, then $\g(S)\cong A_k$.
\label{aklem}\end{lemma}

\begin{proof}
The rank of $S$ is $k$, and the Cartan matrix or Dynkin diagram of the root system is the same as that of $A_k$, thus proving the claim.
\end{proof}

\begin{proof}[of Lemma \ref{annonprim}]
First we consider the cases $n=2,3$ and then show the general result for $n\geq 4$. 
 
\noindent \textbf{Case $\mathbf{n=2}$}

\setcounter{equation}{0}
The weight $\zeta=\mone\aone+\mtwo\atwo$ is a dominant weight if and only if the following inequalities hold:
\begin{eqnarray}
m_2&\leq& 2m_1\\
m_1 &\leq & 2m_2
\end{eqnarray}
For the highest long root $\mu=\alpha_1+\alpha_2$, let $\zeta$ be a highest weight such that $\zeta\succ\mu$ and  $(\zeta;\mu)$ is not primitive.  Then either (a) $\mone=1$ or (b) $\mtwo=1$.

(a)  First let $m_1=1$.  Then $m_2=1$ or $m_2=2$ by the inequalities above.  If $m_2=1$, then $\zeta=\mu$, which is ruled out.  Let $m_2=2$; then $\zeta=\zeta_1=\alpha_1+2\alpha_2$.  We find $\zeta_1-\mu=\alpha_2$ and thus $S=\{\alpha_2\}$ and $\g(S)\cong A_1$.  Relabeling $\alpha_2$ as $\alpha_1$, 
the projection $p:\g\rightarrow \g(S)$ gives $p(\zeta_1)=2\alpha_1$ and $p(\mu)=\alpha_1$.  Then by Lemma \ref{klem} $K_{\zeta_1,\mu}=K_{p(\zeta_1),p(\mu)}=1$.

(b) Similarly, we find that if $\mtwo=1$, then either $\zeta=\mu$, which is ruled out,  or $\zeta=\zeta_2=2\alpha_1+\alpha_2$.  In the second case $\zeta_2-\mu=\alpha_1$, and we conclude that $S=\{\alpha_1\}$ and $\g(S)\cong A_1$.  Since $p(\zeta_2)=2\aone$ and $p(\mu)=\aone$, we again conclude that $K_{\zeta_2,\mu}=1$  by \ref{klem}.

\setcounter{equation}{0}
\noindent \textbf{Case $\mathbf{n=3}$}

A weight $\zeta$ is a dominant weight if and only if the following inequalities hold.
\begin{eqnarray}
m_2&\leq& 2m_1\\
m_1+m_3 &\leq & 2m_2 \\
m_2&\leq & 2m_3
\end{eqnarray}
In this case $\mu=\alpha_1+\alpha_2+\alpha_3$.  Let $\zeta$ be a highest weight such that $(\zeta;\mu)$ is not primitive.  Then $\zeta=m_1\alpha+m_2\alpha_2+m_3\alpha_3$ and one of the following must be true: (a) $m_1=1, (b) \ m_2=1,$ or (c) $m_3=1$.  We will first find all such $\zeta$ and then determine if any give $K_{\zeta,\mu}=1$.

 (a) If $m_1=1$, then $\zeta=\alpha_1+m_2\alpha_2+m_3\alpha_3$ where $m_2$ and $m_3$ satisfy the above inequalities. By (1), $m_2=1$ or $m_2=2$.  If $m_2=1$ we use (2) to conclude that $m_3=1$ and thus $\zeta=\mu$, which is ruled out.  If $m_2=2$ then inequality (2) gives $m_3\leq 3$ resulting in the following dominant weights $\zeta$ such that$(\zeta;\mu)$ is a nonprimitive pair:
\begin{eqnarray*}
\zeta=\zeta_1&=&\alpha_1+2\alpha_2+\alpha_3\\
\zeta=\zeta_2&=&\alpha_1+2\alpha_2+2\alpha_3\\
\zeta=\zeta_3&=&\alpha_1+2\alpha_2+3\alpha_3
\end{eqnarray*}

(b)  Next consider $m_2=1$.  By (2), $m_1=m_3=1$ also and $\zeta=\mu$, which is ruled out.

(c)  Finally, let $m_3=1$.  By (3), $m_2=1$ or $m_2=2$.  Again, if $m_2=1$, then $\zeta=\mu$, which is ruled out.  If $m_2=2$, by (2) $m_1\leq 3$ and thus we are left with two additional dominant weights $\zeta$ such that $(\zeta; \mu)$ is not primitive:
\begin{eqnarray*}
\zeta=\zeta_4&=&2\alpha_1+2\alpha_2+\alpha_3\\
\zeta=\zeta_5&=&3\alpha_1+2\alpha_2+\alpha_3
\end{eqnarray*}
We now determine $K_{\zeta_i,\mu}=1$ in each case.

\begin{enumerate}
\item $\zeta_1-\mu=\alpha_2$.  Thus $S=\{\alpha_2\}$ and $\g(S)\cong A_1$.  Relabeling $\alpha_2$ as $\alpha_1$ we obtain $p(\zeta_1)=2\alpha_1$ and $p(\mu)=\alpha_1$.  We conclude that $K_{\zeta_1,\mu}=1$ by Lemma \ref{klem}.
\item  $\zeta_2-\mu=\alpha_2+\alpha_3$.  Thus $S=\{\alpha_2,\alpha_3\}$ and $\g(S)\cong A_2$.  Relabeling $\{\alpha_2,\alpha_3\}$ as $\{\alpha_1,\alpha_2\}$ we obtain  $p(\zeta_2)=2\alpha_1+2\alpha_2$ and $p(\mu)=\alpha_1+\alpha_2$.  According to Theorem \ref{bzthm} if $K_{\zeta_2,\mu}=K_{p(\zeta_2),p(\mu)}=1$, then $p(\zeta_2)=l\omega_1$ for some positive integer $l$.  We would have $p(\zeta_2)=2\alpha_1+2\alpha_2=l\omega_1=\frac{l}{3}(2\alpha_1+\alpha_2)$.  Each weight is written as the unique sum of simple roots with positive integer coefficients, so this is impossible.  Thus $K_{\zeta_2,\mu}\neq 1$.
\item $\zeta_3-\mu=\alpha_2+2\alpha_3$.  Thus $S=\{\alpha_2,\alpha_3\}$ and $\g(S)\cong A_2$.  Relabeling $\{\alpha_2,\alpha_3\}$ as $\{\alpha_1,\alpha_2\}$ yields $p(\zeta_3)=2\alpha_1+3\alpha_2$ and $p(\mu)=\alpha_1+\alpha_2$.  As above, this satisfies the conditions for $K_{\zeta_3,\mu}=1$ by Theorem \ref{bzthm} if $p(\zeta_3)=2\alpha_1+3\alpha_2=\frac{l}{3}(2\alpha_1+\alpha_2)$ for some positive integer $l$.  However the argument of (2) shows that $l$ does not exist.  Thus $K_{\zeta_3,\mu}\neq 1$.
\item $\zeta_4-\mu=\alpha_1+\alpha_2$.  Thus $S=\{\alpha_1,\alpha_2\}$ and $\g(S)\cong A_2$.  Then $p(\zeta_4)=2\alpha_1+2\alpha_2$ and $p(\mu)=\alpha_1+\alpha_2$ and we have the same conditions as for $\zeta_2$.  Thus the same conclusion holds; $K_{\zeta_4,\mu}\neq 1$.
\item $\zeta_5-\mu=2\alpha_1+\alpha_2$.  Thus $S=\{\alpha_1,\alpha_2\}$ and $\g(S)\cong A_2$.  Then $p(\zeta_5)=3\alpha_1+2\alpha_2$ and $p(\mu)=\alpha_1+\alpha_2$.  As above, this satisfies the conditions for $K_{\zeta_5,\mu}=K_{p(\zeta_5),(\mu)}=1$ by Theorem \ref{bzthm} if $p(\zeta_5)=3\alpha_1+2\alpha_2=\frac{l}{3}(2\alpha_1+\alpha_2)$ for some $l$.  However, again there is no such $l$.  Thus $K_{\zeta_5,\mu}\neq 1$.
\end{enumerate}

Thus for $n=3$, $\zeta=\alpha_1+2\alpha_2+\alpha_3$ is the only highest weight such that $(\zeta;\mu)$ is a nonprimitive pair and $\kzeta=1$.

\noindent \textbf{Case $\mathbf{n\geq 4}$}  

Now we consider the general case for $n\geq 4$.  First we will show that if $\zeta\succ\mu$ and $(\zeta;\mu)$ is not a primitive pair for $\mu=\aone+\cdots+\an$, then $\zeta$ must have one of the following forms:
\begin{eqnarray*}
\zeta=\zeta_1&=&\alpha_1+m_2\alpha_2+m_3\alpha_3+\cdots+m_{n-1}\alpha_{n-1}+\alpha_n\\
\zeta=\zeta_2&=&\alpha_1+m_2\alpha_2+m_3\alpha_3+\cdots+m_{n-1}\alpha_{n-1}+m_n\alpha_n\\
\zeta=\zeta_3&=&m_1\alpha_1+m_2\alpha_2+m_3\alpha_3+\cdots+m_{n-1}\alpha_{n-1}+\alpha_n
\end{eqnarray*}
where $m_i\geq 2$.

Then we will show that $K_{\zeta_i,\mu}\neq 1$ in each of these cases, allowing us to conclude that there are no weights $\zeta$ for $n\geq 4$ such that $(\zeta;\mu)$ is a nonprimitive pair and $K_{\zeta,\mu}=1$.

\setcounter{equation}{0}
For $\zeta=\sum m_i\alpha_i$ a dominant weight, the following inequalities must hold:
\begin{eqnarray}
m_2 &\leq& 2m_1\\
m_{i-1}+m_{i+1}&\leq& 2m_i, \ i=2,\dots,n-1\\
m_{n-1}&\leq& 2m_n
\end{eqnarray}

\begin{sublemma}  Suppose that $m_i=1$ for some $i$, with $i=2,\dots,n-1$.  Then $m_i=1$ for %
all $i=1,\dots,n$.\label{amione}
\end{sublemma}

\begin{proof} Recall that $m_i\geq 1$ for all $i$ by Proposition \ref{atleast1}.  If $\mi=1$ for some $i=2,\dots,n-1$, then by inequality (2), $m_{i-1}+m_{i+1}\leq 2m_i=2$, resulting in $m_{i-1}=m_{i+1}=1$.  By induction on (2), then $m_i=1$ for all $i$.\end{proof}

By the result above, it is clear that $\{\zeta_1,\zeta_2,\zeta_3\}$ are the only dominant weights $\zeta$ different from $\mu$ for which $(\zeta;\mu)$ is not a primitive pair since $\mu=\alpha_1+\cdots+\alpha_n$.

Next we show that $K_{\zeta_i,\mu}\neq 1$ in each case.

\noindent \textbf{Case 1} $\zeta_1=\alpha_1+m_2\alpha_2+\cdots+m_{n-1}\alpha_{n-1}+\alpha_n$, where $m_i\geq 2$ for $2\leq i\leq n-1$.

From the difference $\zeta_1-\mu=(m_2-1)\alpha_2+\cdots +(m_{n-1}-1)\alpha_{n-1}$ we see that  $S=\{\alpha_2,\dots,\alpha_{n-1}\}$ and then $\g(S)\cong A_{n-2}$.  Relabeling $\{\alpha_2,\dots,\alpha_{n-1}\}$ as $\{\alpha_1,\dots,\alpha_{n-2}\}$ we obtain $p(\zeta_1)=m_2\alpha_1+m_3\alpha_2+\cdots+m_{n-1}\alpha_{n-2}$ and $p(\mu)=\alpha_1+\cdots+\alpha_{n-2}$.  We see that  $(p(\zeta_1); p(\mu))$ is now a primitive pair.  By Theorem \ref{bzthm}, $K_{\zeta_1,\mu}=K_{p(\zeta_1),p(\mu)}=1$ if and only if  $p(\zeta_1)=l\omega_1=\frac{l}{n-1}((n-2)\alpha_1+(n-3)\alpha_2+\cdots+\alpha_{n-2})$ for some $l$.  If this is true, then 
\begin{eqnarray*}
m_2&=&\frac{l(n-2)}{n-1}\\
m_3&=&\frac{l(n-3)}{n-1}\\
\vdots\\
m_{n-2}&=&\frac{2l}{n-1}\\
m_{n-1}&=&\frac{l}{n-1}\\
\end{eqnarray*}
or equivalently 
\begin{eqnarray*}
m_{n-2}&=&2m_{n-1}\\
m_{n-3}&=&3m_{n-1}\\
\vdots\\
m_{2}&=&(n-2)m_{n-1}\\
\end{eqnarray*}
Thus we conclude that $m_i=(n-i)m_{n-1}$; i.e. $\{m_i\}_{i=2}^{n-1}$ is strictly decreasing.  However, we have assumed that $m_1=1$, so by inequality (1), $m_2\leq 2$.  Then $m_3<m_2$ means that $m_3\leq 1$ which is a contradiction.  Therefore there is no such dominant weight of the form $\zeta_1$ such that $K_{\zeta_1,\mu}=1$.

\noindent \textbf{Case 2}  $\zeta_2=\alpha_1+m_2\alpha_2+\cdots +m_{n-1}\alpha_{n-1}+m_n\alpha_n$, where  $m_i\geq 2$ for $2\leq i\leq n-1$.

From the difference $\zeta_2-\mu=(m_2-1)\alpha_2+\cdots +(m_{n}-1)\alpha_{n}$ we obtain $S=\{\alpha_2,\dots,\alpha_{n}\}$ and then conclude that $\g(S)\cong A_{n-1}$. Relabeling $\{\alpha_2,\dots,\alpha_n\}$ as $\{\alpha_1,\dots,\alpha_{n-1}\}$ we obtain $p(\zeta_2)=m_2\alpha_1+m_3\alpha_2+\cdots+m_{n}\alpha_{n-1}$ and $p(\mu)=\alpha_1+\cdots+\alpha_{n-1}$.  We see that $(p(\zeta_2); p(\mu))$ is now a primitive pair.  By Theorem \ref{bzthm}, $K_{\zeta_2,\mu}=K_{p(\zeta_2),p(\mu)}=1$ if and only if $p(\zeta_2)=l\omega_1=\frac{l}{n}((n-1)\alpha_1+(n-2)\alpha_2+\cdots+\alpha_{n-1})$ for some $l$.  As above, this requirement allows us to conclude that $\{m_i\}_{i=2}^n$ is a decreasing sequence.  However, we find the same contradiction as above, and therefore $K_{\zeta_2, \mu}\neq 1$.

\noindent \textbf{Case 3}  $\zeta_3=m_1\alpha_1+m_2\alpha_2+\cdots +m_{n-1}\alpha_{n-1}+\alpha_n$, where $m_i\geq 2$ for $2\leq i\leq n-1$. 

From the difference $\zeta_3-\mu=(m_1-1)\alpha_1+\cdots +(m_{n-1}-1)\alpha_{n-1}$ we obtain $S=\{\alpha_1,\dots,\alpha_{n-1}\}$ and then conclude that $\g(S)\cong A_{n-1}$.  Then $p(\zeta_3)=m_1\alpha_1+m_2\alpha_2+\cdots+m_{n-1}\alpha_{n-1}$ and $p(\mu)=\alpha_1+\cdots+\alpha_{n-1}$.  Again, $(p(\zeta_3); p(\mu))$ is now a primitive pair.  By Theorem \ref{bzthm}, $K_{p(\zeta_3),p(\mu)}=K_{\zeta_3,\mu}=1$ if and only if $p(\zeta_3)=l\omega_1=\frac{l}{n}((n-1)\alpha_1+(n-2)\alpha_2+\cdots+\alpha_{n-1})$ for some $l$.  If this is true, then the following set of equalities holds.
\begin{eqnarray*}
m_1&=&\frac{l}{n}(n-1)\\
m_2&=&\frac{l}{n}(n-2)\\
\vdots\\
m_{n-1}&=&\frac{l}{n}\\
\end{eqnarray*}
This is equivalent to
\begin{eqnarray*}
m_1&=&(n-1)m_{n-1}\\
\vdots\\
m_{n-2}&=&2m_{n-1}
\end{eqnarray*}
By inequality (3) $m_{n-1}\leq 2m_n=2$, we conclude that  $m_{n-1}=2$ since $m_{n-1}=1$ implies $m_i=1$ for all $i$ by Sublemma \ref{amione}.  Then $m_i=2(n-i)$ for all $i$ and in particular $m_{n-2}=4$.  However, by (2), $m_{n-2}+m_n \leq 2m_{n-1}$ which implies that $4+1\leq 4$, an obvious contradiction.  Therefore there is no such dominant weight $\zeta_3$ such that $K_{\zeta_3,\mu}=1$, proving our claim.
\end{proof}

\section{Nonprimitive pairs for $B_n, \ n\geq 2$}

In this case, recall that the highest short and long roots are $\mu_1=\alpha_1+\cdots+\alpha_n$ and $\mu_2=\alpha_1+2\alpha_2+2\alpha_3+\cdots+2\alpha_n$.  Since we are restricting to the case that all roots are weights, we only need to find dominant weights $\zeta$ such that $\zeta\succ\mu_2$ and consider nonprimitive pairs $(\zeta;\mu_i)$ for $i=1,2$.  Then we will determine for which of these weights $K_{\zeta,\mu_i}=1$ for $i=1$ or $i=2$.

We recall a result of \cite{decoste} for the $B_n$ case:

\begin{sublemma} Let $\lambda=m_1\alpha_1+\cdots m_n\alpha_n\in\Lambda^+(V)$.  If $m_1=1$, 
then either $m_i=1$ for $1\leq i\leq n$ or $m_i\geq 2$ for $2\leq i \leq n$.\label{eight}\end{sublemma}

\begin{lemma} For $\g$ of type $B_n$, the only highest weights $\zeta$ such that $\zeta\succ\mu_2$, $(\zeta;\mu_i)$ is nonprimitive and  $K_{\zeta,\mu_i}=1$ for $i=1$ or $i=2$ are
\begin{enumerate}
\item $\zeta_1=\alpha_1+2\alpha_2+\mthree\athree$, where $\mthree\geq 3$
\item $\zeta_2=\alpha_1+2\alpha_2+3\alpha_3+3\alpha_4+\cdots+3\alpha_n$, $n\geq 4$  
\item $\zeta_3=2\alpha_1+2\alpha_2+\cdots+2\alpha_n$, $n\geq 2$ 
\end{enumerate}
\label{nonprimb}
\end{lemma}

\begin{proof}
First we consider the case $n=2$ and then show the result for $n\geq 3$.

\noindent \textbf{Case} $\mathbf{n=2}$  

The weight $\zeta=\mone\aone+\mtwo\atwo$ is a dominant weight if and only if the following inequalities hold:
\setcounter{equation}{0}
\begin{eqnarray}
m_2 &\leq & 2m_1\\
m_{1} & \leq & m_2
\end{eqnarray}
In this case $\mu_1=\aone+\atwo$ and $\mu_2=\aone+2\atwo$.  Let $\zeta=\mone\aone+\mtwo\atwo$ be a dominant weight such that $\zeta\succ\mu_2$ and $(\zeta;\mu_i)$ is not a primitive pair for $i=1$ or $i=2$.  Hence $\mone\geq 1$, $\mtwo\geq 2$ and one of the following equalities must hold:  (a) $\mone=1$ or (b) $\mtwo=2$.

(a)  First consider $\mone=1$.  By the inequalities above and the fact that $\zeta\succ\mu_2$, we have $\mtwo=2$.   Thus $\zeta=\mu_2$, which is ruled out.

(b)  Let $\mtwo=2$.  By the inequalities above, $\mone=1$ or $\mone=2$.  If $\mone=1$, we have case (a), so let $\mone=2$.  Then $\zeta=2\aone+2\atwo$ and the pair $(\zeta;\mu_1)$ is primitive while $(\zeta;\mu_2)$ is nonprimitive.  From the difference $\zeta-\mu_2=\aone$, we observe that $S=\{\aone\}$ and $\g(S)\cong A_1$.  The projection $p:\g\rightarrow \g(S)$ gives $p(\zeta)=2\aone$ and $p(\mu_2)=\aone$.  Thus by Lemma \ref{klem}, $K_{\zeta,\mu_2}=1$.  This is $\zeta_3$ above for $n=2$.

\noindent \textbf{Case} $\mathbf{n\geq 3}$

For $\zeta=m_1\alpha_1+\cdots+m_n\alpha_n$ a dominant weight, the following inequalities must hold:
\setcounter{equation}{0}
\begin{eqnarray}
m_2 &\leq & 2m_1\\
m_{i-1}+m_{i+1}&\leq & 2m_i \mbox{ for } i=2,\dots,n-1\\
m_{n-1} & \leq & m_n 
\end{eqnarray}

We find all dominant weights $\zeta=m_1\alpha_1+\cdots+m_n\alpha_n$ such that $\zeta\succ\mu_2$ and  $(\zeta;\mu_i)$ is not a primitive pair for either $i=1$ or $i=2$.  If $\zeta$ is a dominant weight such that $(\zeta;\mu_1)$ is a nonprimitive pair, then $m_i=1$ for some $i$.  If $\zeta$ is a dominant weight such that $(\zeta;\mu_2)$ is a nonprimitive pair, then $m_1=1$ or $m_i=2$ for $i=2,\dots,n$, or both.  Once we have found all such $\zeta$, we will then find $K_{\zeta,\mu_i}$.  

Since $\zeta=\mone\aone+\cdots+\mn\an\succ\mu_2$ we have (*) $\mone\geq 1$ and $\mi\geq 2$ for all $i=2,\dots,n$.

\begin{sublemma} Let $m_i=2$ for some $i\geq 3$, then $m_i=2$ for $i=2,\dots,n$.\end{sublemma}
\begin{proof}
By (2), $m_{i-1}+m_{i+1}\leq 2m_i=4$.  Then if $m_{i-1}>2$ it must be that $m_{i-1}=3$ and $m_{i+1}=1$ which contradicts Sublemma \ref{eight}.  Thus $m_{i-1}=m_{i+1}=2$.  By induction on (2) then $m_i=2$ for $i=2,\dots,n$.  
\end{proof}

If $(\zeta;\mu_1)$ is a nonprimitive pair with $\zeta\succ\mu_2$, then $\mone=1$.  By inequality (1) and the fact that $\mtwo\geq 2$, we conclude that $\mtwo=2$.  Then $\zeta=\aone+2\atwo+\mthree\athree+\cdots+\mn\an$.  If $\mi=2$ for some $i=3,\dots,n$, then by the previous sublemma $\mi=2$ for all $i=3,\dots,n$ which means that  $\zeta=\mu_2$, which we have ruled out.  Hence $\mi\geq 3$ for $i=3,\dots, n$.  We have proved the following:
\begin{enumerate}
\item[(a)]  If $(\zeta; \mu_1)$ is a nonprimitive pair with $\zeta\succ\mu_2$, then $\zeta=\aone+2\atwo+\mthree\athree+\cdots+\mn\an$, where $\mi\geq 3$ for $i\geq 3$.
\end{enumerate}

We turn our attention to $\mu_2$ and find those dominant weights $\zeta$ such that $\zeta\succ \mu_2$ and $(\zeta;\mu_2)$ is not a primitive pair.  We may assume that $\mone\geq 2$ for if $\mone=1$ then $\zeta$ lies in the list (a) by the argument above.

If $\mone\geq 2$, $\zeta\succ\mu_2$ and $(\zeta;\mu_2)$ is not a primitive pair, then $\mi=2$ for some $i\geq 2$.  If $i\geq 3$, then $\zeta=\mu_2$ by the sublemma above, but this is ruled out.  Hence $\mtwo=2$ and from (2) it follows that $2+\mthree\leq \mone+\mthree\leq 2\mtwo=4$.  This implies $\mthree\leq 2$ but $\mthree\geq 2$ since $\zeta\succ\mu_2$.  We conclude that $\mthree=2$ and $\mone=2$.  From the sublemma above we obtain
\begin{enumerate}
\item[(b)]  If $(\zeta; \mu_2)$ is a nonprimitive pair with $\zeta\succ\mu_2$, then either $\zeta$ lies in the list (a) or $\zeta=2\aone+2\atwo+\cdots+2\an$.
\end{enumerate}

Thus we conclude that if $\zeta\succ \mu_2$ is a dominant weight with $(\zeta;\mu_1)$ or $(\zeta;\mu_2)$ nonprimitive, then by (a) and (b) we have two cases to consider.  In the second case we consider only $(\zeta;\mu_2)$.
\begin{eqnarray*}
\zeta &=&\alpha_1+2\alpha_2+m_3\alpha_3+\cdots+m_n\alpha_n,  \ m_i\geq 3 \mbox{ for } i\geq 3\\
\zeta &=& 2\alpha_1+\cdots+2\alpha_n
\end{eqnarray*}

\noindent \textbf{Case 1}  $\zeta=\alpha_1+2\alpha_2+m_3\alpha_3+\cdots+m_n\alpha_n$, where $m_i\geq 3$ for $i\geq 3$.

First consider the nonprimitive pair $(\zeta;\mu_1)$.  From the difference $\zeta-\mu_1=\atwo+(\mthree-1)\athree+\cdots+(\mn-1)\an$, we obtain $S=\{\atwo,\dots,\an\}$ and thus $\g(S)\cong B_{n-1}$ by comparing Dynkin diagrams.  Relabeling $\{\atwo,\dots,\an\}$ as $\{\aone,\dots,\alpha_{n-1}\}$ yields $p(\zeta)=2\aone+\mthree\atwo+\cdots+\mn\alpha_{n-1}$ and $p(\mu_1)=\aone+\cdots+\alpha_{n-1}$.  By Theorem \ref{bzthm}, $K_{p(\zeta),p(\mu_1)}=K_{\zeta,\mu_1}=1$ implies that $p(\mu_1)=\sum_{1\leq i\leq n-1}a_i\omega_i$, where $a_i\in\mathbb{Z}_+$ is  even.  However, $p(\mu_1)=\omega_1$, so $a_1=1$ is not even.  Therefore $K_{\zeta,\mu_1}\neq 1$.

We next consider the nonprimitive pair $(\zeta;\mu_2)$.  From the difference $\zeta-\mu_2=(m_3-2)\alpha_3+\cdots+(m_n-2)\alpha_n$ we see that $S=\{\alpha_3,\dots,\alpha_n\}$.  We need to consider two subcases:  $n=3$ and $n\geq 4$.

\textbf{Subcase $\mathbf{n=3}$}

 Here $\g(S)\cong A_1$.  Relabeling $\{\athree\}$ as $ \{\aone\}$, we obtain $p(\zeta)=\mthree\aone=2\mthree\omega_1$ and $p(\mu_2)=2\aone=4\omega_1$.  By the $A_1$ case of Theorem \ref{bzthm}, $K_{p(\zeta),p(\mu_2)}=1$ if $p(\zeta)=l\omega_1$ and $p(\mu_2)=a\omega_1$ where $l-a\in 2\mathbb{N}$.  Thus $\mu_2$ has multiplicity one if 
$2\mthree-4\in 2\mathbb{N}$, which holds for $\mthree\geq 2$.  However, if $\mthree=2$, then $\zeta=\mu_2$, which is ruled out.  Hence for  $\mthree\geq 3$, $(\zeta;\mu_2)$ is nonprimitive and $K_{\zeta,\mu_2}=1$.  Define $\zeta=\zeta_1$ in this case.

\textbf{Subcase $\mathbf{n\geq 4}$} 

Here $\g(S)\cong B_{n-2}$ since the Dynkin diagrams are the same.  Relabeling $\{\alpha_3,\dots,\alpha_n\}\rightarrow \{\alpha_1,\dots,\alpha_{n-2}\}$ we obtain $p(\zeta)=m_3\alpha_1+m_4\alpha_2+\cdots+m_n\alpha_{n-2}$ and $p(\mu_2)=2\alpha_1+\cdots+2\alpha_{n-2}=2\omega_1$.

We use the $B_n$ case of Theorem \ref{bzthm} to determine if $K_{p(\zeta),p(\mu_2)}=K_{\zeta,\mu_2}=1$, namely if $p(\zeta)=l\omega_1$, $p(\mu_2)=\sum_{1\leq i\leq n-2}a_i\omega_i$, where $a_i\in\mathbb{Z}_+$, even and $(l-1)=\sum_{1\leq i \leq n-3}ia_i+(n-2)a_{n-2}/2$.  

Since $p(\mu_2)=2\omega_1$ it follows that $a_1=2$ and $a_i=0$ for all other $i\neq 1$.  Thus $l=3$ and then for $K_{p(\zeta),p(\mu_2)}=1$ we must have $p(\zeta)=l\omega_1=3\alpha_1+\cdots+3\alpha_{n-2}$; that is, $m_i=3$ for all $i=3,\dots,n$.  Thus in this case, the only dominant weight $\zeta$ such that $K_{\zeta,\mu_2}=1$ is $\zeta=\zeta_2=\alpha_1+2\alpha_2+3\alpha_3+3\alpha_4+\cdots+3\alpha_n$.  

\noindent \textbf{Case 2} $\zeta=\zeta_3=2\alpha_1+\cdots+2\alpha_n$, $n\geq 3$.

From the difference $\zeta-\mu_2=\alpha_1$ we observe that $S=\{\alpha_1\}$ and $\g(S)\cong A_1$.  Then $p(\zeta)=2\alpha_1$ and $p(\mu_2)=\alpha_1$.  By Lemma \ref{klem}, $K_{p(\zeta),p(\mu_2)}=K_{2\alpha_1,\alpha_1}=1$.
\end{proof}

\section{Nonprimitive pairs for $C_n, \ n\geq 3$}

In this case recall that $\mu_1=\aone+2\atwo+\cdots+2\anone+\an$ and $\mu_2=2\aone+2\atwo+\cdots+2\anone+\an$ are the highest short and long roots.

\begin{lemma}  For $\g$ of type $C_n$, there are no highest weights $\zeta$ such that $\zeta\succ\mu_2$, $(\zeta;\mu_i)$ is a nonprimitive pair and $K_{\zeta,\mu_i}=1$ for $i=1$ or $i=2$.\label{nonprimc}\end{lemma}

\begin{proof}

Let us first consider the case $n=3$ and then we will investigate the general case for $n\geq 4$. 

\noindent \textbf{Case $\mathbf{n=3}$}

Here $\mu_1=\aone+2\atwo+\athree$ and $\mu_2=2\aone+2\atwo+\athree$ and we let $\zeta=\mone\aone+\mtwo\atwo+\mthree\athree\succ\mu_2$,  which implies (*) $\mone\geq 2$, $\mtwo\geq 2$ and $\mthree\geq 1$.  If $(\zeta;\mu_i)$ 
is a nonprimitive pair for $i=1,2$, then one of the following must hold:  (a) $\mone=2$, (b) $\mtwo=2$ or (c) $\mthree=1$.  Recall first that $\zeta$ is a dominant weight if and only if the following inequalities hold.
\setcounter{equation}{0}
\begin{eqnarray}
\mtwo&\leq &2\mone\\
\mone+2\mthree &\leq & 2\mtwo\\
\mtwo &\leq & 2\mthree
\end{eqnarray}

We consider each of the cases above to determine nonprimitive pairs.

(a) Suppose that $\mone=2$.  Then by (1) $\mtwo\leq 2\mone=4$ and by (*) $\mtwo\geq 2$, so therefore $\mtwo=2,3$ or $4$.
\begin{enumerate}
\item[(i)] If $\mtwo=2$, then (2)  gives $2+2\mthree=\mone+2\mthree\leq 2\mtwo=4$ and (3)  gives $2=\mtwo\leq 2\mthree$, together yielding $\mthree=1$.  Then $\zeta=\mu_2$, which is ruled out.%
\item[(ii)] If $\mtwo=3$, then by (2) $2+2\mthree=\mone+2\mthree\leq 2\mtwo=6$ and by (3) $3=\mtwo\leq 2\mthree$, giving  $\mthree=2$.  Then $\zeta=2\aone+3\atwo+2\athree$, a candidate.%
\item[(iii)] If $\mtwo=4$, then by (2) $2+2\mthree=\mone+2\mthree\leq 2\mtwo=8$ and by (3) $4=\mtwo\leq 2\mthree$.  We conclude $2\leq \mthree\leq 3$.  Then $\zeta=2\aone+4\atwo+2\athree$ or $\zeta=2\aone+4\atwo+3\athree$, both candidates.
\end{enumerate}

(b)  Suppose that $\mtwo=2$.  By (2) $\mone+2\mthree\leq 2\mtwo=4$ and by (*), $\mone\geq 2$ and $\mthree\geq 1$, forcing the inequalities to be equalities.  We conclude that $\zeta=\mu_2$, which is ruled out.

(c) If $\mthree=1$, then  by (3) $\mtwo\leq 2\mthree=2$ and by (*) $\mtwo\geq 2$.   Thus $\mtwo=2$.  By (2) we see that $\mone\leq 2$ and equality holds by (*).  Again we find that $\zeta=\mu_2$, which is ruled out.

Note that if $\zeta$ is a dominant weight such that $\zeta\succ\mu_2$, then  $(\zeta;\mu_1)$ is a primitive pair by (*) and cases (b) and (c) above.  By the discussion above the only dominant weights $\zeta$ such that $\zeta\succ\mu_2$ and $(\zeta;\mu_2)$ is not a primitive pair are 
\begin{eqnarray*}
\zeta&=&\zeta_1=2\aone+3\atwo+2\athree\\
\zeta&=&\zeta_2=2\aone+4\atwo+2\athree\\
\zeta&=&\zeta_3=2\aone+4\atwo+3\athree
\end{eqnarray*}
We will show for $i=1,2,3$, that $K_{\zeta_i,\mu_2}\neq 1$.  %

From the difference $\zeta_1-\mu_2=\atwo+\athree$, we obtain $S=\{\atwo,\athree\}$ and $\g(S)\cong B_2$.  Then relabeling $\{\atwo,\athree\}$ as $\{\atwo,\aone\}$ yields $p(\zeta_1)=2\aone+3\atwo$ and $p(\mu_2)=\aone+2\atwo$.  By Theorem \ref{bzthm}, for a Lie algebra of type $B_2$, $K_{\zeta_1,\mu_2}=K_{p(\zeta_1),p(\mu_2)}=1$ if and only if $p(\zeta_1)=l\omega_1$ and $p(\mu_2)=a_1\omega_1+a_2\omega_2$ where $a_i\in\mathbb{Z}_+$ are even and $(l-1)=a_1+a_2$.  In our case, $p(\mu_2)=2\omega_2$ and hence $l=3$.  Then for $\mu_2$ to have multiplicity one, $p(\zeta_1)=2\aone+3\atwo=3\omega_1=3(\aone+\atwo)$, which is false.  Therefore $K_{\zeta_1,\mu_2}\neq 1$.

From the difference $\zeta_2-\mu_2=2\atwo+\athree$, we obtain $S=\{\atwo,\athree\}$ and $\g(S)\cong B_2$ as above.  Then relabeling $\{\atwo,\athree\}$ as $\{\atwo,\aone\}$ yields $p(\zeta_2)=2\aone+4\atwo$ and $p(\mu_2)=\aone+2\atwo$.  We are in the same case of Theorem \ref{bzthm} as for $\zeta_1$ where $l=3$.  Then for $\mu_2$ to have multiplicity one, $p(\zeta_2)=2\aone+4\atwo=3\omega_1=3(\aone+\atwo)$, which is again false.  Therefore $K_{\zeta_2,\mu_2}\neq 1$.

From the difference $\zeta_3-\mu_2=2\atwo+2\athree$, we obtain $S=\{\atwo,\athree\}$ and $\g(S)\cong B_2$ as above.  Then relabeling $\{\atwo,\athree\}$ as $\{\atwo,\aone\}$ yields $p(\zeta_3)=3\aone+4\atwo$ and $p(\mu_2)=\aone+2\atwo$.  Again, we are in the same case of \ref{bzthm} as for $\zeta_1$ where $l=3$.  Then for $\mu_3$ to have multiplicity one, $p(\zeta_3)=3\aone+4\atwo=3\omega_1=3(\aone+\atwo)$, which is again false.  Therefore $K_{\zeta_3,\mu_2}\neq 1$ also.

Thus, in the case $n=3$, there are no dominant weights $\zeta$ such that $(\zeta;\mu_2)$ is a nonprimitive pair and  $K_{\zeta,\mu_2}=1$.

\noindent \textbf{Case} $\mathbf{n\geq 4}$ 

If $\zeta=\mone\aone+\cdots+\mn\an$ is a dominant weight, then the following inequalities must hold.
\setcounter{equation}{0}
\begin{eqnarray}
\mtwo &\leq & 2\mone\\
\mimone + \mipone &\leq & 2\mi, \ i=2,\dots,n-2\\
m_{n-2}+2\mn &\leq & 2m_{n-1}\\
m_{n-1} &\leq & 2\mn
\end{eqnarray}
Since we also assume $\zeta\succ\mu_2$, the following inequalities must hold as well: (*) $m_i\geq 2$ for $i=1,\dots,n-1$ and $m_n\geq 1$.  For $\zeta$ such that $(\zeta;\mu_i)$ is not a primitive pair, $i=1,2$, at least one of these inequalities must be an equality.
\begin{sublemma}  If $m_i=2$ for some $i=2,\dots,n-1$, then $m_i=2$ for all $i=1,\dots,n-1$ and $m_n=1$.\label{micnout}
\end{sublemma}

\begin{proof}
Suppose that $m_i=2$ for some $i=2,\dots,n-2$.  Then by inequality (2) and (*), $\mimone=\mipone=2$.  We continue by induction on (2) and find that $m_i=2$ for $i=1,\dots,n-1$.  If $m_{n-1}=2$, then by (3) $m_{n-2}+2m_n\leq 2m_{n-1}=4$.  By (*) we know that  $m_{n-2}\geq 2$ and $m_n\geq 1$, and it follows that $m_{n-2}=2$ and $m_n=1$.  We now apply the first part of the argument to conclude that $m_i=2$ for $i=1,\dots,n-1$.
\end{proof}

If $\zeta$ is a dominant weight such that $\zeta\succ\mu_2$, then we show that $(\zeta;\mu_1)$ is primitive and $(\zeta;\mu_2)$ is nonprimitive only if $\mone=2$.  By Sublemma \ref{micnout} if $m_i=2$ for some $i=2,\dots,n-1$, then $m_i=2$ for $i=1,\dots,n-1$ and $m_n=1$. It follows that $\zeta=\mu_2$, which is ruled out.  Similarly, if $m_n=1$, then by inequality (4) and (*) we obtain $m_{n-1}=2$.  From  Sublemma \ref{micnout} we conclude that $\zeta=\mu_2$, which is ruled out.  Hence $(\zeta;\mu_1)$ is a primitive pair. Thus the only case left to consider is when $m_1=2$. 

If $\zeta\succ\mu_2$ and $(\zeta;\mu_2)$ is nonprimitive, then the previous paragraph shows that $\zeta=2\aone+\mtwo\atwo+\cdots+\mn\an$ where $m_i\geq 3$ for $i=2,\dots,n-1$ and $m_n\geq 2$.  We show that $K_{\zeta,\mu_2}\neq 1$.

From the difference $\zeta-\mu_2=(\mtwo-2)\atwo+\cdots+(m_{n-1}-2)\alpha_{n-1}+(\mn-1)\an$ we find that $S=\{\atwo,\dots,\an\}$ and $\g(S)\cong C_{n-1}$ since $n\geq 4$.  Relabeling $\{\atwo,\dots,\an\}$ as $\{\aone,\dots,\alpha_{n-1}\}$ we find $p(\zeta)=\mtwo\aone+\dots+\mn\alpha_{n-1}$ and $p(\mu_2)=2\aone+\cdots+2\alpha_{n-2}+\alpha_{n-1}$.  Thus $(p(\zeta);p(\mu_2))$ is a primitive pair for a Lie algebra of type $C_{n-1}$.  By Theorem \ref{bzthm}, there are no primitive pairs for $C_n$ such that the weight has multiplicity one.  Therefore $K_{\zeta,\mu_2}\neq 1$.  And finally we conclude that there are no dominant weights $\zeta$ with $(\zeta;\mu_2)$ a nonprimitive pair and $K_{\zeta,\mu_2}=1$.
\end{proof}

\section{Nonprimitive pairs for $D_n$, $n\geq 4$}

In this case recall that $\mu_1=\mu_2=\mu=\alpha_1+2\alpha_2+\cdots+2\alpha_{n-2}+\alpha_{n-1}+\alpha_n$ is the highest short and long root since all roots are the same length.  We find all dominant weights $\zeta$ such that $\zeta\succ\mu$ and $(\zeta;\mu)$ is not a primitive pair.  Then we calculate $K_{\zeta,\mu}$.

\begin{lemma} For $\g$ of type $D_n$, the only highest weights $\zeta$ such that $\zeta\succ\mu$,  $(\zeta;\mu)$ is not a primitive pair and  $K_{\zeta,\mu}=1$ are
\begin{enumerate}
\item for $n=4$ \\
$\zeta_1=\alpha_1+2\alpha_2+2\alpha_3+\alpha_4$\\
 $\zeta_2=\alpha_1+2\alpha_2+\alpha_3+2\alpha_4$\\
 $\zeta_3=2\aone+2\atwo+\athree+\afour$
\item for $n\geq 5$ $\zeta_4=2\alpha_1+\cdots+2\alpha_{n-2}+\alpha_{n-1}+\alpha_n$
\end{enumerate}
\label{nonprimd}
\end{lemma}

\begin{proof}
First we consider the case $n=4$ and then show the general result.

\noindent \textbf{Case $\mathbf{n=4}$} 

In this case $\mu=\alpha_1+2\alpha_2+\alpha_3+\alpha_4$.  The weight $\zeta=m_1\alpha_1+m_2\alpha_2+m_3\alpha_3+m_4\alpha_4$ is a dominant weight if and only if the following inequalities hold.
\setcounter{equation}{0}
\begin{eqnarray}
m_2 &\leq & 2m_1\\
m_1+m_3+m_4 &\leq & 2m_2\\
m_2 & \leq & 2m_3\\
m_2 & \leq & 2m_4
\end{eqnarray}
For $\zeta$ such that $\zeta\succ\mu$ and  $(\zeta;\mu)$ is not a primitive pair we have  (*) $\mone\geq 1$, $\mtwo\geq 2$, $\mthree\geq 1$ and $\mfour\geq 1$.  Hence at least one of the following must hold: (a) $m_1=1$, (b) $m_2=2$, (c)  $m_3=1$ or  (d) $m_4=1$.  We consider each case.

(a) If $m_1=1$, then $m_2=2$ by (1) and (*) and $\mthree+\mfour\leq 3$ by (2).  Hence we have one of the following three cases: $m_3=m_4=1$, $m_3=2$ and $m_4=1$ or $m_3=1$ and $m_4=2$.  The corresponding weights are respectively $\zeta=\mu$, which is ruled out, $\zeta=\zeta_1=\alpha_1+2\alpha_2+2\alpha_3+\alpha_4$ and $\zeta=\zeta_2=\alpha_1+2\alpha_2+\alpha_3+2\alpha_4$, respectively.

(b) If $m_2=2$ then $m_1=1$ or $m_1=2$ since $\mone\leq 2$ by (2).  If $m_1=1$, we have the previous case, so let $m_1=2$.  Then by (2), $m_3=m_4=1$ and $\zeta=\zeta_3=2\alpha_1+2\alpha_2+\alpha_3+\alpha_4$.

(c) and (d) If either $\mthree=1$ or $\mfour=1$, then $\mtwo=2$ by (3) or (4) respectively and (*).  Then we have $\zeta$ as in the previous cases.

Thus for $n=4$ we have 3 dominant weights $\zeta$ such that $(\zeta;\mu)$ is not a primitive pair:
\begin{eqnarray*}
\zeta=\zeta_1&=&\alpha_1+2\alpha_2+2\alpha_3+\alpha_4\\
\zeta=\zeta_2&=&\alpha_1+2\alpha_2+\alpha_3+2\alpha_4\\
\zeta=\zeta_3&=&2\alpha_1+2\alpha_2+\alpha_3+\alpha_4
\end{eqnarray*}
In each case the difference $\zeta_i-\mu$ yields $S=\{\alpha_i\}$ and $g(S)\cong A_1$.  Then $p(\xi_1)=2\alpha_1$ and $p(\mu)=\alpha_1$.  By Lemma \ref{klem} we conclude that $K_{\zeta_i,\mu}=1$ in each case. 





\noindent \textbf{Case $\mathbf{n\geq 5}$}  

For $\zeta=m_1\alpha_1+\cdots+m_n\alpha_n$ a dominant weight in this case the following inequalities must hold:
\setcounter{equation}{0}
\begin{eqnarray}
m_2 &\leq& 2m_1\\
m_{i-1}+m_{i+1}&\leq& 2m_i, \ i=2,\dots,n-3\\
m_{n-3}+m_{n-1}+m_n &\leq & 2m_{n-2}\\
m_{n-2}& \leq & 2m_{n-1}\\
m_{n-2} &\leq & 2m_n
\end{eqnarray}
We show that the following are the only dominant weights $\zeta=m_1\alpha_1+\cdots+m_n\alpha_n$ such that $\zeta\succ\mu$ and  $(\zeta;\mu)$ is a nonprimitive pair.  Then we consider $K_{\zeta,\mu}$ in each case.
\begin{eqnarray*}
\zeta=\zeta_1&=&2\alpha_1+\cdots+2\alpha_{n-2}+\alpha_{n-1}+\alpha_n\\
\zeta=\zeta_2&=&\alpha_1+2\alpha_2+3\alpha_3+m_4\alpha_4+\cdots+m_n\alpha_n, \mbox{ with } m_i\geq 3,  \ i=4,\dots,n-2, \\
& & \mbox{ and } m_i\geq 2, \ i=n-1,n
\end{eqnarray*}

 Note that  a dominant weight $\zeta$ such that $\zeta\succ\mu$ and $(\zeta;\mu)$ is a nonprimitive pair must satisfy (*) $m_i\geq 1$ for $i=1,n-1,n$ and  $m_i\geq 2$ for $i=2,\dots,n-2$ with at least one of  these inequalities an equality for $i=1,\dots,n$.  We consider each case below.

\begin{sublemma}  Suppose that $m_i=2$ for some $i$, $i=3,\dots,n-2$, then $m_i=2$ for all $i=2,\dots, n-2$.
\label{midn}
\end{sublemma}

\begin{proof}  Suppose that $m_i=2$ for some $i=3,\dots,n-3$.  Then by (2), $m_{i-1}+m_{i+1}\leq 4$, but $m_{i-1}\geq 2$ and $m_{i+1}\geq 2$ by (*) and therefore $m_{i-1}=m_{i+1}=2$.  By induction on (2), $m_i=2$ for $i=2,\dots,n-2$.  If $m_{n-2}=2$, then by (*) and (3)  $2+m_{n-1}+\mn\leq m_{n-3}+m_{n-1}+\mn\leq 2m_{n-2}=4$.  This implies  that $m_{n-1}=\mn=1$ and $m_{n-3}=2$.  By the previous case $\mi=2$ for $i=2,\dots,n-2$.
\end{proof}


\noindent \textbf{Case} $m_{n-1}=1$ or $m_n=1$.    

If $m_{n-1}=1$, then by (4) and (*), $m_{n-2}=2$ and by (3) and (*), we have $m_{n-3}=2$ and $m_n=1$.  Similarly, if $m_n=1$, then we conclude that $m_{n-1}=1$ and $m_{n-3}=2$.  It now follows from Sublemma \ref{midn} that if $m_{n-1}=1$ or $\mn=1$ then $m_{n-1}=\mn=1$ and $m_i=2$ for $i=2,\dots,n-2$. 

From (2) we see that $\mone\leq 2$ since $\mone+2=\mone+\mthree\leq 2\mtwo=4$.  Hence either  $m_1=1$, in which case $\zeta=\mu$, which is ruled out, or $m_1=2$ for which $\zeta=\zeta_1=2\alpha_1+2\alpha_2+\cdots+2\alpha_{n-2}+\alpha_{n-1}+\alpha_n$.  Then $\zeta_1$ is a dominant weight such that $(\zeta_1;\mu)$ is not a primitive pair.

\noindent \textbf{Case} $m_i=2$ for some $i=3,\dots,n-2$.  

Then by Sublemma \ref{midn}, $m_i=2$ for all $i=2,\dots,n-2$.  By (3) $m_{n-1}=m_n=1$ and we have the same result as in the previous case.

\noindent \textbf{Case} $m_2=2$.

If $m_3=2$ also, then we have the previous case, so  we may assume that $m_i\geq 3$ for $i=3,\dots,n-2$.  We show that there exists a dominant weight $\zeta_2$ distinct from the previous weight $\zeta_1$ such that $(\zeta_2;\mu)$ is nonprimitive.  By (4) $3\leq m_{n-2}\leq 2m_{n-1}$ and by (5) $3\leq m_{n-2}\leq 2m_n$, giving $m_{n-1}\geq 2$ and $m_n\geq 2$.  By (2) we have $m_1+m_3\leq 2\mtwo=4$, but we also know that $m_3\geq 3$ and therefore $m_1=1$ and $m_3=3$.  Thus $\zeta_2=\alpha_1+2\alpha_2+3\alpha_3+m_4\alpha_4+\cdots+m_n\alpha_n$ is a dominant weight such that $(\zeta_2;\mu)$ is not a primitive pair if $\mi\geq 3$ for $i=4,\dots,n-2$ and $\mi\geq 2$ for $i=n-1,n$.

\noindent \textbf{Case} $\mone=1$.

 By (1) and (*) $\mtwo=2$ and we have the same result as in the previous case.

Thus the following are the only 2 dominant weights $\zeta$ such that $(\zeta;\mu)$ is not a primitive pair:
\begin{eqnarray*}
\zeta=\zeta_1&=&2\alpha_1+\cdots+2\alpha_{n-2}+\alpha_{n-1}+\alpha_n\\
\zeta=\zeta_2&=&\alpha_1+2\alpha_2+3\alpha_3+m_4\alpha_4+\cdots+m_n\alpha_n, \mbox{ with } m_i\geq 3,  \ i=4,\dots,n-2, \\
& & \mbox{ and } m_i\geq 2, \ i=n-1,n
\end{eqnarray*}

We find $K_{\zeta_i,\mu}$ in each case.

From the difference $\zeta_1-\mu=\alpha_1$, we see that $S=\{\alpha_1\}$ and $\g(S)\cong A_1$.  Then $p(\zeta_1)=2\alpha_1$ and $p(\mu)=\alpha_1$ and by Lemma \ref{klem}, $K_{\zeta_1,\mu}=1$.  Note:  this is $\zeta_4$ of Lemma \ref{nonprimd}, which is the same as $\zeta_3$ for $n=4$.

Next we consider $\zeta_2$, first for the case $n=5$.  In this case, $\zeta_2=\alpha_1+2\alpha_2+3\alpha_3+m_4\alpha_4+m_5\alpha_5$, where $\mfour\geq 2$ and $\mfive\geq 2$,  and $\mu=\alpha_1+2\alpha_2+2\alpha_3+\alpha_4+\alpha_5$.  The inequality (3) gives $2+m_4+m_5\leq 2\mthree=6$ and we conclude then that the only possible values for $m_4$ and $m_5$ are $m_4=m_5=2$.  In this case $\zeta_2=\alpha_1+2\alpha_2+3\alpha_3+2\alpha_4+2\alpha_5$.

From the difference $\zeta_2-\mu=\alpha_3+\alpha_4+\alpha_5$ we see that $S=\{\alpha_3,\alpha_4,\alpha_5\}$ and through the relabeling $\{\alpha_3,\alpha_4,\alpha_5\}$ as $\{\alpha_2,\alpha_1,\alpha_3\}$ we find $\g(S)\cong A_3$.  This yields $p(\zeta_2)=2\alpha_1+3\alpha_2+2\alpha_3$ and $p(\mu)=\alpha_1+2\alpha_2+\alpha_3$.  We see that $(p(\zeta_2);p(\mu))$ is a primitive pair and we can now apply Theorem \ref{bzthm} to determine if $K_{\zeta_2,\mu}=K_{p(\zeta_2),p(\mu)}=1$.  According to this result, if $K_{p(\zeta_2),p(\mu)}=1$, then $p(\zeta_2)=l\omega_1$.  Hence $p(\zeta_2)=2\alpha_1+3\alpha_2+2\alpha_3=\frac{l}{4}(3\alpha_1+2\alpha_2+\alpha_3)$ for some positive integer $l$.  However, there is no such $l$ and therefore $K_{\zeta_2,\mu}\neq 1$.

Next, consider $n\geq 6$.  From the difference $\zeta_2-\mu=\alpha_3+(m_4-2)\alpha_4+\cdots+(m_n-1)\alpha_n$ we observe that $S=\{\alpha_3,\dots,\alpha_n\}$ and since $n\geq 6$, the Dynkin diagram gives $\g(S)\cong D_{n-2}$.  Relabeling $\{\alpha_3,\dots,\alpha_n\}$ as $\{\alpha_1,\dots,\alpha_{n-2}\}$ we obtain $p(\zeta_2)=3\alpha_1+m_4\alpha_2+\cdots+m_n\alpha_{n-2}$ and $p(\mu)=2\alpha_1+\cdots+2\alpha_{n-4}+\alpha_{n-3}+\alpha_{n-2}$. We see that $(p(\zeta_2);p(\mu))$ is a primitive pair in $D_{n-2}$, but by Theorem \ref{bzthm} we also observe that there are no primitive pairs $(\zeta;\mu)$ in $D_{n-2}$ such that $\kzeta=1$, therefore $K_{\zeta_2,\mu}\neq 1$. \end{proof}

\section{Nonprimitive pairs for $E_n$}
\subsection{Nonprimitive pairs for $E_6$}

In this case, recall that $\mu_1=\mu_2=\mu=\aone+2\atwo+2\athree+3\afour+2\afive+\asix$ is the highest root.  

\begin{lemma}  For $\g$ of type $E_6$ there are no highest weights $\zeta$ with $\zeta\succ\mu$ and $(\zeta;\mu)$ a nonprimitive pair such that $\kzeta=1$.
\label{nonprimesix}\end{lemma}
\begin{proof}

Recall that any dominant weight $\zeta=\mone\aone+\cdots+\msix\asix$ must satisfy the following:
\setcounter{equation}{0}
\begin{eqnarray}
\mthree &\leq & 2\mone\\
\mfour &\leq & 2\mtwo\\
\mone + \mfour & \leq & 2 \mthree\\
\mtwo + \mthree+\mfive &\leq & 2\mfour\\
\mfour + \msix &\leq & 2\mfive\\
\mfive &\leq & 2\msix
\end{eqnarray}

Note that $\zeta\succ \mu$ yields the following inequalities for $\zeta=\mone\aone+\cdots+\msix\asix$: (*) $\mi\geq 1$ for $i=1,6$, $\mi\geq 2$ for $i=2,3,5$ and $\mfour\geq 3$.  In addition, if $\zeta$ is a dominant weight such that $(\zeta;\mu)$ is not a primitive pair, at least one of the inequalities will be an equality for some $i$.  We consider each of these cases.

\noindent \textbf{Case 1}  Suppose that $\mone=1$.  By (1) and (*) $\mthree=2$ and  by (3) $1+\mfour\leq 4$, thus by (*), $\mfour=3$.  By (4), $\mtwo+2+\mfive\leq 2\mfour=6$ which with (*) implies that $\mtwo=\mfive=2$.  Lastly, by (5), $\msix=1$.  Thus if $\mone=1$, then $\zeta=\mu$, which is ruled out.

\noindent \textbf{Case 2}  Next, suppose that $\mtwo=2$.  By (2) $\mfour\leq 4$, so $\mfour=3$ or $\mfour=4$ by (*).  If $\mfour=3$, then by (4) and (*) we conclude that $\mtwo=\mthree=\mfive=2$, which implies that $\mone=1$ by (3).  Again, $\msix=1$ by (5) and we find that $\zeta=\mu$, which is ruled out.

Let $\mfour=4$.  By (3) we have $\mone+4\leq 2\mthree$ and by (5) $4+\msix\leq 2\mfive$ which give $\mthree\geq 3$ and $\mfive \geq 3$.  Inequality (4) yields $\mthree+\mfive\leq 6$ and thus we conclude that $\mthree=\mfive=3$.  Note that $\mone=2$ follows from (1) and (3).  Also, $\msix=2$ by (5) and (6).  Thus $\zeta=\zeta_1=2\aone+2\atwo+3\athree+4\afour+3\afive+2\asix$ is a dominant weight such that $(\zeta_1;\mu)$ is not a primitive pair.

\noindent \textbf{Case 3}  Suppose that $\mthree=2$.  By (3) and (*) $\mone=1$ and $\mfour=3$.  Therefore we are in the first case considered and $\zeta=\mu$, which we have ruled out.

\noindent \textbf{Case 4}  Suppose that $\mfour=3$.  Then by (2), $\mtwo\geq 2$, but by (4) and (*), $\mtwo=\mthree=\mfive=2$.  As seen above, if $\mtwo=2$, then $\zeta=\mu$ or $\zeta=\zeta_1$, however the identity $\mfour=3$ yields only $\zeta=\mu$, which is ruled out.

\noindent \textbf{Case 5}  Suppose that $\mfive=2$.  By (5) $\mfour+\msix\leq 4$, and thus $\mfour=3$ and $\msix=1$ by (*).  We are now in the previous case.

\noindent \textbf{Case 6}  Suppose that $\msix=1$.  Then by (6) and (*) $\mfive=2$ and we are in the previous case.

From the above, we conclude that  $\zeta=2\aone+2\atwo+3\athree+4\afour+3\afive+2\asix$ is the only dominant weight such that $(\zeta;\mu)$ is a nonprimitive pair.  We will show that $K_{\zeta,\mu}\neq 1$.

From the difference $\zeta-\mu=\aone+\athree+\afour+\afive+\asix$ we observe that $S=\{\aone,\athree,\afour,\afive,\asix\}$ and by considering the Dynkin diagram of $E_6$, we see that $\g(S)\cong A_5$.  Relabeling $\{\aone,\athree,\afour,\afive,\asix\}$ as $\{\aone,\atwo,\athree,\afour,\afive\}$, we obtain $p(\zeta)=2\aone+3\atwo+4\athree+3\afour+2\afive$ and $p(\mu)=\aone+2\atwo+3\athree+2\afour+\afive$.  We notice that $(p(\zeta);p(\mu))$ is a primitive pair and therefore we apply Theorem \ref{bzthm} to determine if $K_{p(\zeta),p(\mu)}=1$.  In the $A_5$ case, if $K_{p(\zeta),p(\mu)}=1$, then $p(\zeta)=l\omega_1=\frac{l}{6}(5\aone+4\atwo+3\athree+2\afour+\afive)$ for some positive integer $l$.  Clearly, there is no such $l$ and therefore $K_{\zeta,\mu}\neq 1$. \end{proof}

\subsection{Nonprimitive pairs for $E_7$}
In this case, recall that $\mu_1=\mu_2=\mu=2\aone+2\atwo+3\athree+4\afour+3\afive+2\asix+\aseven$ is the highest root.  
\begin{lemma}  For $\g$ of type $E_7$ there are no highest weights $\zeta$ with $\zeta\succ\mu$ and $(\zeta;\mu)$ a nonprimitive pair such that $\kzeta=1$.
\label{nonprimeseven}\end{lemma}
\begin{proof}

Recall that any dominant weight $\zeta=\mone\aone+\cdots+\mseven\aseven$ must satisfy the following:
\setcounter{equation}{0}
\begin{eqnarray}
\mthree &\leq & 2\mone\\
\mfour &\leq & 2\mtwo\\
\mone + \mfour & \leq & 2 \mthree\\
\mtwo + \mthree+\mfive &\leq & 2\mfour\\
\mfour + \msix &\leq & 2\mfive\\
\mfive +\mseven &\leq & 2\msix\\
\msix &\leq & 2\mseven
\end{eqnarray}
Note that $\zeta\succ \mu$ yields the following inequalities for $\zeta=\mone\aone+\cdots+\mseven\aseven$:  (*) $\mi\geq 2$ for $i=1,2,6$, $\mi\geq 3$ for $i=3,5$, $\mfour\geq 4$ and $\mseven\geq 1$.  In addition, if $\zeta$ is a dominant weight such that $(\zeta;\mu)$ is not a primitive pair, at least one of the inequalities will be an equality for some $i$.  We consider each of these cases.

We show that the inequalities above can  be satisfied for $\zeta\succ\mu$, $(\zeta;\mu)$ a nonprimitive pair, only if $\zeta=\zeta_1=2\aone+3\atwo+4\athree+6\afour+5\afive+4\asix+2\aseven$ or $\zeta=\zeta_2=2\aone+3\atwo+4\athree+6\afour+5\afive+4\asix+3\aseven$.  We then show that $K_{\zeta_i,\mu}\neq 1$ for $i=1,2$.

\noindent \textbf{Case 1}  Suppose that $\mone=2$.  By (1) and (*) $\mthree=3$ or $\mthree=4$.

First suppose that $\mthree=3$.  By (3) $2+\mfour\leq 6$ and combined with (*), we conclude $\mfour=4$.  Inequality (4) then gives $\mtwo+3+\mfive\leq 8$ and since $\mtwo\geq 2$ and $\mfive \geq 3$ by (*), these are in fact equalities.  From (5) and (*), we see that $\msix=2$ which means  that $\mseven=1$ by (6).  Thus if $\mone=2$ and $\mthree=3$, then $\zeta=\mu$, which is ruled out.

Now let $\mthree=4$.  By inequality (3)  $2+\mfour\leq 8$, giving $\mfour\leq 6$.  By (*) $\mfour\geq 4$.  We consider each case $m=4,5,6$ individually.
\begin{enumerate}
\item[$\mathbf{\mfour=4}$] Inequality (4) and (*) give $2+4+3\leq \mtwo+\mthree+\mfive\leq 2\mfour=8$, an obvious contradiction.  Therefore, $\mfour\neq 4$.
\item[$\mathbf{\mfour=5}$]  By (2) $5=\mfour\leq 2\mtwo$ and thus $\mtwo\geq 3$.  Also, by (5) and (*), $5+2\leq \mfour+\msix\leq 2\mfive$, so $\mfive\geq 4$.  We get a contradiction by (4) since then $3+4+4\leq \mtwo+\mthree+\mfive\leq 2\mfour=10$.  Therefore, $\mfour\neq 5$.
\item[$\mathbf{\mfour=6}$]  We consider the inequalities in the case $\mfour=6$ where $\mthree=4$ and $\mone=2$.  We show that $\zeta=\zeta_1$ or $\zeta=\zeta_2$ as listed above.
\begin{enumerate}
\item $\mtwo\geq 3$\\
This follows from (2) since $\mfour=6$.
\item $\mfive=5$\\
From (a) and (4) we have $7+\mfive\leq \mtwo+\mthree+\mfive\leq 2\mfour=12$.  Hence $\mfive\leq 5$.  If $\mfive\leq 4$, then by (5) we have $6+\msix=\mfour+\msix\leq 2\mfive\leq 8$, which implies $\msix\leq 2$.  By (*) $\msix\geq 2$, and hence $\msix=2$ and $\mfive=4$.  By (6) we obtain $4+\mseven=\mfive+\mseven\leq 2\msix=4$, which is impossible.  Hence $\mfive=5$.
\item $\mtwo=3$\\
From (4) and (b) we have $\mtwo+9=\mtwo+\mthree+\mfive\leq 2\mfour=12$.  Hence $\mtwo\leq 3$ and equality holds by (a).
\item $\msix=4$\\
From (5) and (b) we have $6+\msix=\mfour+\msix\leq 2\mfive=10$, which implies $\msix\leq 4$.  If $\msix\leq 3$, then by (6) and (b) we have $5+\mseven=\mfive+\mseven\leq 2\msix\leq 6$.  This implies $\mseven\leq 1$, and equality holds by (*).  This implies that $\msix=3$, but by (7) we have $\msix\leq 2\mseven=2$.  This contradiction shows that $\msix=4$.
\item $\mseven=2$ or 3\\
By (7) and (d) we have $4=\msix\leq 2\mseven$, which implies that $\mseven\geq 2$.   By (6), (b) and (d) we have $5+\mseven=\mfive+\mseven\leq 2\msix=8$, which implies that $\mseven\leq 3$.
\end{enumerate}
\end{enumerate}

We then have 2 dominant weights $\zeta_i$ such that $(\zeta_i;\mu)$ are not primitive pairs:
$\zeta_1= 2\aone+3\atwo+4\athree+6\afour+5\afive+4\asix+2\aseven$ and 
$\zeta_2= 2\aone+3\atwo+4\athree+6\afour+5\afive+4\asix+3\aseven$.

\noindent \textbf{Case 2}  Suppose that $\mtwo=2$.  Then by (2) and (*) $4\leq \mfour\leq 2\mtwo=4$, so $\mfour=4$.  By (4) and (*), $2+3+3\leq \mtwo+\mthree+\mfive\leq 2\mfour=8$, which results in $\mthree=\mfive=3$.  Then by (3) and (*) $2+4\leq \mone+\mfour\leq 2\mthree=6$ and therefore $\mone=2$.  We are now in Case 1 with $\mthree=3$, but this was ruled out.

\noindent \textbf{Case 3}  Suppose that $\mthree=3$.  By (3) and (*) $2+4\leq \mone+\mfour\leq 2\mthree=6$, and hence $\mone=2$ and $\mfour=4$.  We are again in Case 1 with $\mthree=3$, which was ruled out.

\noindent \textbf{Case 4}  Suppose that $\mfour=4$.  By (4) and (*) $2+3+3\leq \mtwo+\mthree+\mfive\leq 2\mfour=8$, yielding $\mtwo=2$, $\mthree=3$, and $\mfive=3$.  We are now in Case 2, which was ruled out.

\noindent \textbf{Case 5}  Let $\mfive=3$.  By (5) and (*) $4+2\leq \mfour+\msix\leq 2\mfive=6$, and hence $\mfour=4$ and $\msix=2$.   We are now in Case 4, which was ruled out.

\noindent  \textbf{Case 6}  Let $\msix=2$.  Then by (6) and (*) $3+1\leq \mfive+\mseven \leq 2\msix=4$, yielding the equalities $\mfive=3$ and $\mseven=1$.  We are in Case 5, which was ruled out.

\noindent \textbf{Case 7}  Let $\mseven=1$.  By (7) and (*) $2\leq \msix\leq 2\mseven=2$ which implies that $\msix=2$.  Again, we fall into the previous case which was ruled out.

Thus, the only two dominant weights $\zeta$ such that $\zeta\succ\mu$ and $(\zeta;\mu)$ is not a primitive pair are $\zeta_1$ and $\zeta_2$ as in Case 1.  
We now consider $K_{\zeta_i,\mu}$ in each case.

From the difference $\zeta_1-\mu=\atwo+\athree+2\afour+2\afive+2\asix+\aseven$ we observe that $S=\{\atwo,\athree,\afour,\afive,\asix,\aseven\}$ and by considering the Dynkin diagram of $E_7$, we see that $\g(S)\cong D_6$.  Then relabeling $\{\atwo,\athree,\afour,\afive,\asix,\aseven\}$ as $ \{\asix,\afive,\afour,\athree,\atwo,\aone\}$, we obtain $p(\zeta_1)=2\aone+4\atwo+5\athree+6\afour+4\afive+3\asix$ and $p(\mu)=\aone+2\atwo+3\athree+4\afour+3\afive+2\asix$.  Then $(p(\zeta_1);p(\mu))$ is a primitive pair for $D_6$.  By Theorem \ref{bzthm}, for type $D_6$, there are no primitive pairs such that the dimension of the weight space is one, so therefore $K_{\zeta_1,\mu}\neq 1$.

The case of $\zeta_2$ is similar to the previous one.  From the difference $\zeta_2-\mu=\atwo+\athree+2\afour+2\afive+2\asix+2\aseven$ we observe that again $S=\{\atwo,\athree,\afour,\afive,\asix,\aseven\}$ and then $\g(S)\cong D_6$.  Then with the same relabeling $\{\atwo,\athree,\afour,\afive,\asix,\aseven\}$ as $\{\asix,\afive,\afour,\athree,\atwo,\aone\}$, we obtain $p(\zeta_2)=3\aone+4\atwo+5\athree+6\afour+4\afive+3\asix$ and $p(\mu)=\aone+2\atwo+3\athree+4\afour+3\afive+2\asix$.  Then $(p(\zeta_2);p(\mu))$ is a primitive pair for $D_6$, and by Theorem \ref{bzthm}, we conclude that $K_{\zeta_2,\mu}\neq 1$.  
\end{proof}

\subsection{Nonprimitive pairs for $E_8$}
In this case, recall that $\mu_1=\mu_2=\mu=2\aone+3\atwo+4\athree+6\afour+5\afive+4\asix+3\aseven+2\aeight$.  

\begin{lemma}  For $\g$ of type $E_8$ there are no highest weights $\zeta\succ \mu$ such that $(\zeta;\mu)$ is a nonprimitive pair with $\kzeta=1$.
\label{nonprimeeight}\end{lemma}

\begin{proof}

Recall that any dominant weight $\zeta=\mone\aone+\cdots+\meight\aeight$ must satisfy the following:
\setcounter{equation}{0}
\begin{eqnarray}
\mthree &\leq & 2\mone\\
\mfour &\leq & 2\mtwo\\
\mone + \mfour & \leq & 2 \mthree\\
\mtwo + \mthree+\mfive &\leq & 2\mfour\\
\mfour + \msix &\leq & 2\mfive\\
\mfive +\mseven &\leq & 2\msix\\
\msix +\meight &\leq & 2\mseven\\
\mseven &\leq & 2\meight
\end{eqnarray}

Note that $\zeta\succ \mu$ yields the following inequalities for $\zeta=\mone\aone+\cdots+\meight\aeight$: (*) $\mi\geq 2$ for $i=1,8$, $\mi\geq 3$ for $i=2,7$, $\mi\geq 4$ for $i=3,6$, $\mfour\geq 6$ and $\mfive\geq 5$.  In addition, if $\zeta$ is a dominant weight such that $(\zeta;\mu)$ is not a primitive pair, at least one of the inequalities will be an equality for some $i$.  We consider each of these cases.

\noindent \textbf{Case 1}  Let $\mone=2$.  By (1) $\mthree\leq 2\mone=4$ and by (*) $\mthree\geq 4$, so we conclude that $\mthree=4$. By (3) and (*) it then  follows that $2+6\leq m_1+m_4\leq 2\mthree=8$, resulting in  $\mfour=6$.  We find $\mtwo=3$ and $\mfive=5$ by (*) and (4).  From (5) and (*) we obtain  $\msix=4$.  From (6) and (*) we see that 
 $\mseven=3$.  By (7) and (*) we obtain 
$\meight=2$.  Hence $\mone=2$ implies that $\zeta=\mu$, which is ruled out.

\noindent \textbf{Case 2}  Let $\mtwo=3$.  By (2) and (*) 
$\mfour=6$.  By (4) and (*) 
$\mthree=4$ and $\mfive=5$.  By (3) and (*) we have 
$\mone=2$.  By the first case $\zeta=\mu$, which is ruled out.

\noindent \textbf{Case 3}  Let $\mthree=4$.  Then (3) and (*) 
force $\mone=2$ and $\mfour=6$. By Case 1 we have $\zeta=\mu$, which is ruled out.

\noindent \textbf{Case 4}  Let $\mfour=6$. Then  (4) and (*) 
yield $\mtwo=3$, $\mthree=4$ and $\mfive=5$.  By Case 2, we have $\zeta=\mu$, which is ruled out.

\noindent \textbf{Case 5}  Let $\mfive=5$.  By (5) and (*) 
$\mfour=6$ and $\msix=4$.  By the previous case $\zeta=\mu$, which is ruled out.

\noindent \textbf{Case 6}  Let $\msix=4$.  By (6) and (*) 
$\mfive=5$ and $\mseven=3$. We are now in Case 5.

\noindent \textbf{Case 7}  Let $\mseven=3$.  By (7) and (*) 
$\msix=4$ and $\meight=2$.  We are now in Case 6.

\noindent \textbf{Case 8}  Let $\zeta\succ\mu$ be a dominant weight with $(\zeta;\mu)$ a nonprimitive pair and $\meight=2$.  Since the previous seven cases have been ruled  out we may assume that $\mone\geq 3$, $\mtwo\geq 4$, $\mthree\geq 5$, $\mfour\geq 7$, $\mfive\geq 6$, $\msix\geq 5$ and $\mseven\geq 4$.  It follows that $(p(\zeta);p(\mu))$ is a primitive pair, where $S=\{\aone,\dots,\aseven\}$.  Hence $\g(S)\cong E_7$ by an inspection of the Dynkin diagram.  However by Theorem \ref{bzthm} no primitive pairs $(p(\zeta);p(\mu))$ with $K_{p(\zeta),p(\mu)}=1$  exist for $E_7$.

\vspace{0.1in}

\noindent \textbf{Remark}  With further work one can show that the only dominant weight $\zeta$ with $\zeta\succ\mu$, $\meight=2$ and $(\zeta;\mu)$ a nonprimitive pair is $\zeta=4\aone+5\atwo+7\athree+10\afour+8\afive+6\asix+4\aseven+2\aeight$.
\end{proof}

\section{Nonprimitive pairs for $F_4$}

Recall that in this case the highest short and long roots are $\mu_1=\aone+2\atwo+3\athree+2\afour$ and $\mu_2=2\aone+3\atwo+4\athree+2\afour$.  

\begin{lemma}  For $\g$ of type $F_4$ there are no highest weights $\zeta$ such that $\zeta\succ\mu_2$, $(\zeta;\mu_i)$ is a nonprimitive pair and $K_{\zeta,\mu_i}=1$ for $i=1,2$.\label{nonprimf}\end{lemma}

\begin{proof}

Recall that a weight $\zeta$ is a dominant weight if and only if the following inequalities hold:
\setcounter{equation}{0}
\begin{eqnarray}
\mtwo&\leq&2\mone\\
\mone+\mthree&\leq&2\mtwo\\
2\mtwo+\mfour&\leq&2\mthree\\
\mthree&\leq&2\mfour
\end{eqnarray}

Since $\zeta\succ\mu_2$ the following inequalities must hold: (*) $\mone\geq 2$, $\mtwo\geq 3$, $\mthree\geq 4$ and $\mfour \geq 2$.  In addition if $\zeta$ is a dominant weight such that the pair is not primitive, at least one of the inequalities must be an equality.  We consider each case.

\noindent \textbf{Case 1}  Let $\mone=2$.  Then by (1) and (*) $3\leq \mtwo\leq 4$.  
\begin{enumerate}
\item[(a)]  If $\mtwo=3$ , then by (2) and (*) $m_3=4$.  Then by (3) and (*) $\mfour=2$.  Hence $\zeta=\mu_2$, which is ruled out.
\item[(b)]  If $\mtwo=4$, then by (2)  $\mthree\leq 6$.  By (3) $\mthree\geq 5$.  Therefore $\mthree=5$ or $\mthree=6$.  If $\mthree=5$ then by (3) $\mfour=2$, resulting in a contradiction in inequality(4): $5=\mthree\leq 2\mfour=4$.  Hence $\mthree=6$.  Inequality (3) implies that  $\mfour\leq 4$.  By (4) $6=\mthree\leq 2\mfour$, which implies that $\mfour=3$ or $\mfour=4$.
\end{enumerate}
Thus in Case 1 there are two dominant weights $\zeta_i$ such that $\zeta_i\succ\mu_2$ and $(\zeta_i;\mu_2)$ is not a primitive pair, namely $\zeta_1=2\aone+4\atwo+6\athree+3\afour$ and $\zeta_2=2\aone+4\atwo+6\athree+4\afour$.

\noindent \textbf{Case 2}  Let $\mtwo=3$.  Then by (2) and (*)  $\mone=2$ and $\mthree=4$.  We are in (a) of Case 1, which is ruled out.

\noindent \textbf{Case 3}  Let $\mthree=4$.  Then by (3) and (*)  $\mtwo=3$ and $\mfour=2$.  By (2) and (*) $\mone=2$ and $\zeta=\mu_2$, which is ruled out.

\noindent \textbf{Case 4}  Let $\mfour=2$.  By (4) and (*) we have $4\leq \mthree\leq 2\mfour=4$ and thus $\mthree=4$ as in the previous case, forcing $\zeta=\mu_2$, which is ruled out.

Thus we have found two dominant weights $\zeta_i$ such that $(\zeta_i;\mu_2)$ is a nonprimitive pair: 
\begin{eqnarray*}
\zeta_1&=&2\aone+4\atwo+6\athree+3\afour\\
\zeta_2&=&2\aone+4\atwo+6\athree+4\afour
\end{eqnarray*}
Next we show that $K_{\zeta_i,\mu_2}\neq 1$ in each case.  Note that by Theorem \ref{bzthm} $K_{\zeta_i,\mu_1}\neq 1$ since $(\zeta_i;\mu_1)$ is a primitive pair for $i=1,2$.

From the difference $\zeta_1-\mu_2=\atwo+2\athree+\afour$ we observe that $S=\{\atwo,\athree,\afour\}$ and by comparing Dynkin diagrams that $\g(S)\cong C_3$.  Relabeling $\{\atwo,\athree,\afour\}$ as $\{\athree,\atwo,\aone\}$ yields $p(\zeta_1)=3\aone+6\atwo+4\athree$ and $p(\mu_2)=2\aone+4\atwo+3\athree$.  Now $(p(\zeta_1);p(\mu_2))$ is a primitive pair in $C_3$ and by Theorem \ref{bzthm} there are no primitive pairs for $C_3$ such that the weight space has dimension one. Therefore $K_{\zeta_1,\mu_2}\neq 1$.

From the difference $\zeta_2-\mu_2=\atwo+2\athree+2\afour$ we observe that again $S=\{\atwo,\athree,\afour\}$ and $\g(S)\cong C_3$.  Relabeling $\{\atwo,\athree,\afour\}$ as $\{\athree,\atwo,\aone\}$ yields $p(\zeta_2)=4\aone+6\atwo+4\athree$ and $(\mu_2)=2\aone+4\atwo+3\athree$.  Now $(p(\zeta_2);p(\mu_2))$ is a primitive pair in $C_3$ and similarly by Theorem \ref{bzthm} we conclude that $K_{\zeta_1,\mu_2}\neq 1$.

Thus there are no dominant weights $\zeta$ such that $(\zeta;\mu)$ is a primitive pair and $\kzeta=1$ for $F_4$.
\end{proof}

\section{Nonprimitive pairs for $G_2$}

Recall that the highest short and long roots in this case are $\mu_1=2\aone+\atwo$ and $\mu_2=3\aone+2\atwo$.  

\begin{lemma}  For $\g$ of type $G_2$, the only highest weight $\zeta$ such that $\zeta\succ\mu_2$ and $(\zeta;\mu_i)$ is nonprimitive for $i=1$ or $2$ is $\zeta=4\aone+2\atwo$.  In this case $(\zeta;\mu_1)$ is  primitive and $(\zeta;\mu_2)$ is nonprimitive with $K_{\zeta,\mu_2}=1$.\label{nonprimg}\end{lemma}
\begin{proof}
For $\zeta=\mone\aone+\mtwo\atwo$ to be a dominant weight, the following must hold:
\setcounter{equation}{0}
\begin{eqnarray}
3\mtwo&\leq&2\mone\\
\mone&\leq&2\mtwo
\end{eqnarray}

We are looking for dominant weights $\zeta\succ\mu_2$ such that $(\zeta;\mu_i)$ is a nonprimitive pair for $i=1$ or $2$.  The condition $\zeta\succ\mu_2$ implies that (*) $\mone\geq 3$ and $\mtwo\geq 2$, and hence $(\zeta;\mu_1)$ will be a primitive pair.  We then restrict our discussion to finding dominant weights $\zeta\succ\mu_2$ such that $(\zeta;\mu_2)$ is not a primitive pair.  

For $\zeta$ a dominant weight $\zeta=\mone\aone+\mtwo\atwo$ such that $(\zeta;\mu_2)$ is not a primitive pair, either $\mone=3$ or $\mtwo=2$.  If $\mone=3$, then $\mtwo\geq 2$ by (*) and $\mtwo\leq 2$ by (1).  Hence $\mtwo=2$ and $\zeta=\mu_2$, which is ruled out.  Thus we consider $\mtwo=2$.  In this case inequalities (1) and (2) yield either $\mone=3$ or $\mone=4$ .  In the first case we again have $\zeta=\mu_2$, which is ruled out, but in the second, we have $\zeta=4\alpha_1+2\alpha_2$, and $(\zeta;\mu_2)$ is not a  primitive pair.  We consider $K_{\zeta,\mu_2}$.

From the difference $\zeta-\mu_2=\alpha_1$ we see that $S=\{\alpha_1\}$ and then $\g(S)\cong A_1$.  Then $p(\zeta)=4\aone$ and $p(\mu_2)=3\aone$.  By Theorem \ref{bzthm} for a Lie algebra of type $A_1$, $K_{\zeta,\mu_2}=K_{p(\zeta),p(\mu_2)}=1$ if $p(\zeta)=l\omega_1$, $p(\mu_2)=a_1\omega_1$ where $a_1\in\mathbb{Z}_+$ and $(l-a_1)\in 2\mathbb{N}$.  In this case, $p(\mu_2)=3\aone=6\omega_1$, so $l$ must be an integer such that $(l-6)\in 2\mathbb{N}$ and $p(\zeta)=l\omega_1=\frac{l}{2}\alpha_1$.  Clearly $l=8$ satisfies these conditions and therefore $K_{p(\zeta),p(\mu_2)}=K_{\zeta,\mu_2}=1$.  
\end{proof}

\end{document}